\documentclass[12pt,final]{amsart}

\usepackage[margin=1in]{geometry}

\usepackage{amssymb}

\newtheorem{lemma}{Lemma}[section]

\newtheorem{theorem}[lemma]{Theorem}

\theoremstyle{definition}

\theoremstyle{remark}
\newtheorem{remark}[lemma]{Remark}

\newcommand{\N}{\ensuremath{{\mathbb N}}}

\newcommand{\C}{\ensuremath{\mathbb C}}

\newcommand{\Po}{\ensuremath{\mathcal {P}}}

\newcommand{\nd}{\ensuremath{\mathrm {End}}}
\newcommand{\tr}{\ensuremath{\mathrm {tr}}}

\author{Tatyana Barron}
\address{T. Barron (corresponding author), Department of Mathematics, 
University of Western Ontario, 
London, Ontario N6A 5B7, Canada }
\email{tatyana.barron@uwo.ca}

\author{Noah Wheatley}
\address{ }
\email{nwheatle@uwo.ca}

\thanks{N.W. was a research student of T.B. Research is supported in part 
by the Natural Sciences and Engineering Research Council of Canada}

\title[]{Entanglement and products}

\begin{document}
\sloppy

\maketitle

\noindent {\bf Abstract.} 
We address a general question whether geometry of submanifolds of an integral compact K\"ahler manifold is characterized by invariants that come from analysis. In geometric quantization, we have an integral compact K\"ahler manifold $M$ and a holomorphic line bundle $L$ on this manifold. There is a known procedure how to associate a sequence of mixed states $(\rho_N)$, $N=1,2,3,...$, to a submanifold $\Lambda$ of $M$. Do analytic properties of this sequence reflect the geometry of $\Lambda$ ?  

In this paper, we consider the case when $M$ is a product of two integral compact K\"ahler manifolds. We show that, when $\Lambda$ is a product submanifold of $M$, then the entanglement of formation of $\rho_N$ is zero for all sufficiently large $N$.

\

\noindent {\bf 2020 MSC}:  53D50, 53B50, 81S10, 81P42. 

\

\noindent {\bf Keywords}: Hilbert spaces, mixed states, entanglement of formation, K\"ahler quantization, submanifolds, line bundles. 

\

\section{Introduction}

Many questions in analysis on powers of holomorphic line bundles on K\"ahler manifolds 
are related to quantization. In this context quantization is
understood as a map from the algebra of smooth functions on the K\"ahler manifold
to an algebra of operators on a Hilbert space. The unit vectors in this Hilbert space
are pure states of the quantum system.

In quantum information theory, there is a rich collection of techniques that are available
to study entanglement on tensor products of  Hilbert spaces. 
Berezin-Toeplitz quantization produces sequences of Hilbert spaces that come from the geometry of a symplectic manifold. 
Applying the existing methods of quantifying entanglement to these specific Hilbert spaces 
may lead to insights about quantum states.   
In previous work \cite{barron:17} T. Barron and T. Pollock calculated the entanglement
entropy for specific pure states of geometric nature and concluded that these states
have maximal entanglement entropy. This is meaningful, because,  entropy of entanglement, considered as a resource,
quantifies the amount of information that can be transmitted \cite{nielsen:00}
(with an appropriate use of the quantum system). Entanglement entropy is a quantity that is well defined for pure states. 
For mixed states, typically other entanglement measures are used (e.g. entanglement of formation).
Results in the literature include lower bounds of entanglement of formation \cite{carlen:12}, \cite{carlen:16}.  

In the general context of exploring the interplay between entanglement and symplectic geometry, it is possible to provide a  characterization of entanglement using symplectic structure of the space of quantum states (both for pure and mixed states)  \cite{ozm:14, saw:11}. 

In geometric quantization, there is a recurring general theme of exploring the correspondence between submanifolds (or, more generally, subsets) 
of the phase space, and the corresponding states. This procedure is especially illuminating for Lagrangian states. The paper \cite{barron:17} was an attempt 
to add an entanglement entropy calculation to that process. This aligns with the general effort to bring techniques of quantum information theory over to the setting of K\"ahler quantization.  In particular, \cite{lefloch:18} provides bounds for fidelity of Lagrangian states.  
Entropy, specifically with ties to Berezin-Toeplitz operators, was recently used in \cite{charles:18}, \cite{perez:18}.

In this paper, we investigate the correspondence between geometry of the classical phase space and entropy in the following aspect: we prove that if the classical phase space $M$ is a product submanifold, i.e. $M=M_1\times M_2$, 
then certain quantum states, attached to $M$, are not entangled. For the accurate statement see Theorem \ref{mainth} and see remarks after the Theorem for discussion.  These mixed states are associated with a submanifold of a K\"ahler manifold and previously appeared  in \cite{perez:18}.

\section{Background}

First, we collect some necessary definitions and facts. For quantum information theory, possible references are  \cite{be:17}, \cite{nielsen:00}, \cite{witten:20}. Linear algebra material in this section is standard. 

Let $V_1$ be a $d_1$-dimensional complex vector space and 
$V_2$ be a $d_2$-dimensional complex vector space. W.l.o.g. $d_1\le d_2$. 
 A vector in $V_1\otimes V_2$ is called {\it decomposable} if it can be written as $a\otimes b$, with $a\in V_1$, $b\in V_2$.  
The {\it partial traces} are the linear maps ${\mathrm {Tr}}_1$ and ${\mathrm {Tr}}_2$ defined by 
$$
{\mathrm {Tr}}_1:{\mathrm {End}}(V_1\otimes V_2)\to {\mathrm {End}}(V_2)
$$ 
$$
{\mathrm {Tr}}_1(A\otimes B)={\mathrm{Tr}}(A)B 
$$
$$
{\mathrm {Tr}}_1(\sum_j A_j\otimes B_j)=\sum_j{\mathrm{Tr}}_1(A_j\otimes B_j)
$$
$$
{\mathrm {Tr}}_2:{\mathrm {End}}(V_1\otimes V_2)\to {\mathrm {End}}(V_1)
$$ 
$$
{\mathrm {Tr}}_2(A\otimes B)={\mathrm{Tr}}(B)A 
$$
$$
{\mathrm {Tr}}_2(\sum_j A_j\otimes B_j)=\sum_j{\mathrm{Tr}}_2(A_j\otimes B_j)
$$
for $A,A_j\in{\mathrm {End}}(V_1)$, $B, B_j\in{\mathrm {End}}(V_2)$, and finite sums $\sum\limits_j$. 

Suppose $V_1$ and $V_2$ are equipped with Hermitian inner products. Thus,  each of $V_1$, $V_2$, $V_1\otimes V_2$
is a finite-dimensional Hilbert space. The inner product on $V_1\otimes V_2$ is defined on decomposable vectors by 
$$
\langle a\otimes b, c\otimes d\rangle = \langle a,c\rangle\langle b,d\rangle
$$
 and extended by linearity. Since $V_1$, $V_2$ are finite-dimensional, every linear operator in ${\mathrm {End}}(W)$, where $W$ is one of $V_1$, $V_2$,  
 $V_1\otimes V_2$, is bounded. 
 For a nonzero vector $a$ in $W$,  we will denote by $\Po _a$ the orthogonal projection from $W$ onto the one-dimensional subspace of $W$  spanned by $a$. So, 
 if $||a||=1$, then 
 $$
 \Po _a:   w\mapsto \langle w,a\rangle a, \ w\in W. 
 $$
Let $v$ be a unit vector in $V_1\otimes V_2$.  The {\it entanglement entropy} $\nu(v)$ of $v$ is 
$$
\nu(v)=-\sum_{j=1}^{d_1}\lambda_j\ln \lambda_j
$$
where $\lambda_1$,...,$\lambda_{d_1}$ are the eigenvalues of ${\mathrm {Tr}}_2(\Po _v)$,
with the convention
$0\ln 0=0$.

We can represent $v$ as  
$v=\sum_{j=1}^{d_1}\alpha_je_j\otimes f_j$, where $\{e_i\}$ is an orthonormal basis in $V_1$, $\{ f_j\}$ is an orthonormal set in $V_2$, and
$0\le \alpha_1\le ...\le \alpha_{d_1}$ are real numbers ({\it Schmidt decomposition} of $v$, which is conceptually a restatement of singular value decomposition).   
The numbers $\alpha_j$ are uniquely determined by $v$. The choice of $\{e_i\}$ and $\{f_j\}$ is not unique. The eigenvalues $\lambda_1$,...,$\lambda_{d_1}$ of ${\mathrm {Tr}}_2(\Po _v)$ 
are the numbers $(\alpha_1)^2$,..., $(\alpha_{d_1})^2$. 

The operator $\Po _v$ is Hermitian, positive semidefinite, of trace $1$, and so is  ${\mathrm {Tr}}_2(\Po _v)$. A Hermitian positive semidefinite operator $A\in 
{\mathrm{End}}(W)$ of trace $1$ is a rank $1$ projection (i.e. $A=\Po _v$ for some unit vector $v\in W$) if and only if $A^2=A$. Suppose $v_1,v_2$ are unit vectors in $V_1,V_2$, respectively, and $v=v_1\otimes v_2$. 
Then $\Po _v=\Po _{v_1}\otimes \Po _{v_2}$.  Conversely, if $A\in {\mathrm{End}}(V_1)$ and $B \in {\mathrm{End}}(V_2)$  
are such that $C=A\otimes B$ is Hermitian, positive semidefinite, of trace $1$, and $C^2=C$, then $C=\Po _v$, where $v=a\otimes b$, for some $a\in V_1$, 
$b\in V_2$, such that $||a||=||b||=1$. 

The value of $\nu (v)$ is a real number in the interval $[0,\ln d_1]$. It is zero if and only
if $v$ is decomposable, also if and only if
the eigenvalues of ${\mathrm {Tr}}_2(\Po _v)$ are $1$, $0$,..., $0$.
The value of $\nu(v)$ is $\ln d_1$ if and only if the eigenvalues of ${\mathrm {Tr}}_2(\Po _v)$ are all equal:
$\lambda_1=...=\lambda_{d_1}=\frac{1}{d_1}$. The number $\nu(v)$ characterizes ``how far" the vector $v$
is from being decomposable. This is not meant  
 with respect to the distance induced by the inner product, or another metric. Rather, it quantifies
``how non-decomposable'' the vector $v$ is by characterizing how close the eigenvalues of ${\mathrm {Tr}}_2(\Po _v)$ 
are to being equidistributed. They are all equal if and only if the entanglement entropy is maximal.

Now, suppose $C\in {\mathrm{End}}(V_1\otimes V_2)$ is Hermitian, positive semidefinite, of trace $1$. The {\it entanglement of formation} of $C$ is 
defined as follows \cite{bennett:99, be:17}. If $C$ is a rank $1$ projection (i.e. $C=\Po _v$ for some unit vector $v\in V_1\otimes V_2$), then 
$\tilde{\nu}(C)=\nu(v)$. Otherwise, 
$$
\tilde{\nu}(C) =\inf \sum_j p_j \nu(v_j)
$$
where the sum $\underset{j}{\sum}$ is finite, $p_j$ are positive real numbers whose sum is $1$, the vectors $v_j$ are unit vectors in $V_1\otimes V_2$ such that 
$$
C=\sum_j p_j\Po _{v_j}
$$
and the infimum is taken over all such finite sums.   Note: if $C$ is a rank $1$ projection, then we also can define $\tilde{\nu}(C)$ this way, and it will be consistent with the definition given above.  

Let's clarify. Suppose $C\in {\mathrm{End}}(V_1\otimes V_2)$ is Hermitian, positive semidefinite, of trace $1$. The linear operator $C$ is $\Po _v$ for some unit vector $v\in V_1\otimes V_2$ if and only if $C^2=C$. Suppose $C=\Po _v$. Then $\nu(v)=0$ if and only if $v$ is decomposable. 

Now, suppose  $C\in {\mathrm{End}}(V_1\otimes V_2)$ is Hermitian, positive semidefinite, of trace $1$, but $C^2\ne C$. In this case, it does not make sense to talk about entanglement entropy, and, instead, we use the entanglement of formation $\tilde {\nu}(C)$. 

As a trivial observation, every linear operator $C\in {\mathrm{End}}(V_1\otimes V_2)$ can be written as a finite linear combination $C=\sum_jA_j\otimes B_j$, where $A_j\in{\mathrm{End}}(V_1)$ and $B_j\in{\mathrm{End}}(V_2)$ for all $j$. 

If $C=\sum\limits_j p_jA_j\otimes B_j$ (a finite sum), where, for each $j$, $p_j>0$ and   $A_j\in{\mathrm{End}}(V_1)$ and $B_j\in{\mathrm{End}}(V_2)$ are Hermitian, positive semidefinite,  of trace $1$, then $\tilde {\nu}(C)=0$ \cite{be:17}. 

In physics terminology, a unit vector $v\in V_1\otimes V_2$, or the associated projector $\Po _v\in {\mathrm{End}}(V_1\otimes V_2)$, represents a {\it pure} state in 
the Hilbert space $V_1\otimes V_2$. If $v$ is decomposable, this state is {\it separable (non-entangled)}. Otherwise, this pure state is {\it entangled}.

A  linear operator $C\in {\mathrm{End}}(V_1\otimes V_2)$, which is  Hermitian, positive semidefinite, of trace $1$, represents 
a {\it mixed} state. This state is {\it separable (non-entangled)}, if it is a finite linear combination of separable pure states with positive coefficients 
whose sum is $1$ - i.e., such $C$ is separable, if it can be written in the form 
$$
C=\sum_jp_j \Po _{a_j}\otimes \Po _{b_j}
$$
where the sum is finite, $p_j$ are positive numbers such that $\sum_j p_j=1$, the vectors $a_j\in V_1$ and $b_j\in V_2$ for all $j$. (Note: equivalently, a separable mixed state can be defined as a finite linear combination $\sum_jp_jA_j\otimes B_j$, where all $p_j>0$ and $A_j\in{\mathrm{End}}(V_1)$, 
$B_j\in{\mathrm{End}}(V_2)$ are Hermitian, positive semidefinite, of trace $1$). 
). Otherwise, if the mixed state $C$ is not separable, it is called {\it entangled}. 

We finish this section with a couple of remarks. 

\begin{remark} It is easy to provide an example of a state with a nonzero entanglement of formation: e.g. $\Po _v$, where $v$ is a nondecomposable unit vector is such a mixed state.  Other examples of states with nonzero entanglement of formation can be found in \cite{wootters:98},  
\cite[16.9]{be:17}.  
\end{remark} 
\begin{remark}
From the perspective of applications to the Hilbert spaces that appear in Berezin-Toeplitz quantization, one might try to do the following. Let $V$ be a finite-dimensional Hilbert space. Also pick an inner product on $V^*$. Let $v\in V$ be a vector of norm $1$. There is an isomorphism 
$$
\iota: {\mathrm{End}}(V)\to V^*\otimes V   . 
$$
Calculate the entanglement entropy of the vector $\iota(\Po _v)$. However, the answer will always be zero, and so an attempt to use entropy this way is not fruitful. 
\end{remark} 

\section{Setting and main theorem}

For $j=1,2$ let $(L_j,h_j)$ be a holomorphic hermitian line bundle on a compact connected
K\"ahler manifold  $(M_j,\omega_j)$ of complex dimension $n_j\ge 1$
such that the curvature of the Chern connection on $L_j$ 
is $-i\omega_j$. W.l.o.g. $n_1\le n_2$. 
Denote by 
$dV_j$ the measure on $M_j$ associated to the volume form $\dfrac{w_j^n}{n_j!}$
and set $dV=dV_1dV_2$. 
Let $N$ be a positive integer. 
Recall that the holomorphic line bundle $L_1^N\boxtimes L_2^N\to M_1\times M_2$ is defined by 
$L_1^N\boxtimes L_2^N=\pi_1^*(L_1^N)\otimes \pi_2^*(L_2^N)$, 
where $\pi_1:M_1\times M_2\to M_1$ and  $\pi_2:M_1\times M_2\to M_2$ are the 
projections onto the first and onto the second factor respectively. We have: 
$\pi_1^*(L_1^N)\otimes \pi_2^*(L_2^N)\simeq (\pi_1^*L_1)^N\otimes (\pi_2^*L_2)^N$ (because the tensor product commutes with pullbacks). 

Denote by $h_1^{(N)}$, $h_2^{(N)}$, $h^{(N)}$ the  hermitian metrics on the line bundles
$L_1^N$, $L_2^N$, $L_1^N\boxtimes L_2^N$ respectively, induced by
the hermitian metrics $h_1$ and $h_2$.

To simplify notations, on a specific complex vector space, we will denote by $\langle .,.\rangle$ the appropriate hermitian inner product (if there is no ambiguity). In particular, 
if $z\in M_j$ and $u,v\in (L_j)_z$, then
$$
\langle u,v\rangle =h_j(u,v)
$$
and if $s,\tau\in H^0(M_j,L_j^N)$, then 
$$
\langle s,\tau\rangle =\int_{M_j}h^{(N)}_j(s(z),\tau(z))dV_j(z). 
$$
The holomorphic line bundles $L_1$, $L_2$ are ample.  For sufficiently large $N$, 
$$
\dim H^0(M_1,L_1^N)=C_1N^{n_1}+O(N^{n_1-1})
$$
and 
$$
\dim H^0(M_2,L_2^N)=C_2N^{n_2}+O(N^{n_2-1})
$$
where $C_j=\int_{M_j}\frac{c_1(L_j)}{n_j!}$ for $j=1,2$. 

We will denote 
$$
d_N=\dim H^0(M_1\times M_2,L_1^N \boxtimes L_2^N)
$$
and
$$
d_N^{(j)}=\dim H^0(M_j,L_j^N)
$$
for $j=1,2$. 
So, for sufficiently large $N$, we have 
$d_N^{(1)}\le d_N^{(2)}$. 
\begin{lemma}
  For sufficiently large $N$ (specifically, $N$ such that $L_1^N$ and $L_2^N$ are very ample)  
  $$
  H^0(M_1\times M_2,L_1^N \boxtimes L_2^N)\cong
  H^0(M_1,L_1^N)\otimes H^0(M_2,L_2^N)
  $$
  (isomorphism of Hilbert spaces). 
\end{lemma}
\proof
The obvious $\C$-linear map 
$$
H^0(M_1,L_1^N)\otimes H^0(M_2,L_2^N)\to H^0(M_1\times M_2,L_1^N \boxtimes L_2^N)
$$ 
defined on decomposable elements by
$$
s_1\otimes s_2\mapsto \pi_1^*s_1\otimes \pi_2^*s_2
$$
preserves the inner product. To show that for sufficiently large $N$   
this map is surjective, we verify that 
these two complex vector spaces have the same dimension: 
$$
\dim H^0(M_1\times M_2,L_1^N \boxtimes L_2^N)=\int_{M_1\times M_2} Td(M_1\times M_2)ch (L_1^N \boxtimes L_2^N)=
$$
$$
\int_{M_1}\pi_1^* Td(M_1) \ \pi_1^*ch(L_1^N)\int_{M_2}\pi_2^* Td(M_2) \ \pi_2^*ch(L_2^N)=
\dim H^0(M_1,L_1^N) \ \dim H^0(M_2,L_2^N)=
$$
$$
\dim (H^0(M_1,L_1^N)\otimes H^0(M_2,L_2^N)).
$$
$\Box$

For each $j\in\{ 1,2\}$, choose and fix a point $p_j\in M_j$ and a unit vector $\xi_j\in (L_j)_{p_j}$. Let $N\in\N$. The map  
$$
s\mapsto h^{(N)}_{j}(s(p_j),\xi_j^{\otimes N})
$$
is a continuous linear functional on $H^0(M_j,L_j^N)$ (linearity follows from the linearity of the inner product, and this functional is bounded, because, if $\{ e_k^{(j)} \}$ is a basis 
in  $H^0(M_j,L_j^N)$, then 
$$
s=\sum_{k=1}^{d_N^{(j)}} c_ke_k^{(j)}
$$ 
for some $c_k\in\C$ and 
$$
|h^{(N)}_{j}(s(p_j),\xi_j^{\otimes N})|^2=|h^{(N)}_{j}(\sum_kc_ke_k^{(j)}(p_j),\xi_j^{\otimes N})|^2=
|\sum_k c_kh^{(N)}_{j}(e_k^{(j)}(p_j),\xi_j^{\otimes N})|^2
$$
which, by the Cauchy-Schwarz inequality, does not exceed $C||s||^2$, where 
$$
||s||=\sqrt{\sum_k|c_k|^2}
$$ and 
$$
C=\sqrt { \sum_k|h^{(N)}_{j}(e_k^{(j)}(p_j),\xi_j^{\otimes N})|^2} \ ).
$$
Therefore, there is 
a unique element of $H^0(M_j,L_j^N)$, which we will denote by $\Theta ^{(N)}_{p_j,\xi_j}$, with the property
\begin{equation}
  \label{thetaprop12}
\langle s,\Theta ^{(N)}_{p_j,\xi_j}\rangle=h^{(N)}_j(s(p_j),\xi_j^{\otimes N})
\end{equation}
for each $s\in H^0(M_j,L_j^N)$. We observe:
 for $\lambda\in\C$ such that $|\lambda|=1$ 
\begin{equation}
  \label{eqlin}
\Theta ^{(N)}_{p_j,\lambda\xi_j}=\lambda^N\Theta ^{(N)}_{p_j,\xi_j}.
\end{equation}
Let $p=(p_1,p_2)\in M_1\times M_2$. Let $\xi_1\in(\pi_1^*L_1)_p$, $\xi_2\in(\pi_2^*L_2)_p$ be vectors of norm $1$, and $\xi=\xi_1\otimes \xi_2$. As above, the continuous linear functional
on $H^0(M_1\times M_2,L_1^N \boxtimes L_2^N)$ defined by 
$$
s\mapsto h^{(N)}(s,\xi_1^{\otimes N}\otimes \xi_2^{\otimes N})
$$
determines the unique element $\Theta ^{(N)}_{p;\xi}$ of $H^0(M_1\times M_2,L_1^N \boxtimes L_2^N)$ with the property 
\begin{equation}
  \label{thetaprop}
  \langle s,\Theta ^{(N)}_{p;\xi}\rangle=h^{(N)}(s(p),\xi_1^{\otimes N}\otimes\xi_2^{\otimes N})
\end{equation}
for all $s\in H^0(M_1\times M_2,L_1^N \boxtimes L_2^N)$.

Let $\Lambda$ be a connected smooth submanifold of $M_1\times M_2$. For each $N\in\N$ let 
$$
{\mathcal{R}}_N:H^0(M_1\times M_2,L_1^N \boxtimes L_2^N)\to L^2(\Lambda,(L_1^N \boxtimes L_2^N)\Bigr | _{\Lambda})
$$
$$
s\mapsto s\Bigr |_{\Lambda}
$$
be the restriction operator. Denote by $V_N$ the image of $H^0(M_1\times M_2,L_1^N \boxtimes L_2^N)$ under ${\mathcal{R}}_N$. Let 
$P_N^*$ be the Hilbert space adjoint of the operator 
$$
P_N: H^0(M_1\times M_2,L_1^N \boxtimes L_2^N)\to V_N
$$
$$
s\mapsto {\mathcal{R}}_N(s).
$$
Let    
$$
\rho_N:=\frac{1}{\tr (P_N^*P_N)}P_N^*P_N\in \nd (H^0(M_1\times M_2,L_1^N \boxtimes L_2^N)). 
$$
Now that we fixed all notations, we can state the main theorem.  
\begin{theorem}
  \label{mainth}
$(i)$ Suppose $\Lambda$ is a point $\{p\}$, $p=(p_1,p_2)\in M_1\times M_2$.
For $j=1,2$, let $\xi_j$ be a vector in  $(\pi_j^* L_j)_{p}$ of norm $1$. 
Then $\rho_N$ is the (rank $1$)
orthogonal projection onto the linear subspace spanned by $\Theta_{p;\xi}^{(N)}$.   
For sufficiently large $N$,
the entanglement entropy of the pure state  $\rho_N$
is equal to zero. 

\noindent $(ii)$ Suppose $\Lambda=\Lambda_1\times \Lambda_2$ is a connected smooth submanifold of $M_1\times M_2$ such that 
$\Lambda_1$ is a submanifold of $M_1$, $\Lambda_2$ is a submanifold of $M_2$, and  
$\dim (\Lambda)>0$. 

Then  
for sufficiently large $N$  
$$
\tilde{\nu}(\rho_N)=0.
$$
\end{theorem}
\begin{remark}
  Because of (\ref{eqlin}), the statement $(i)$ in the Theorem does not depend on the choice of $\xi_1$ and $\xi_2$ (meaning that the one-dimensional subspace 
  is the same, regardless of the choice of $\xi_1$, $\xi_2$).   
\end{remark}
\proof
To prove $(i)$, we observe that  
the operator
$$
P^*_N:\C\cong (L_1^N \boxtimes L_2^N)\Bigr | _{p}\to H^0(M_1\times M_2,L_1^N \boxtimes L_2^N)
$$
is the linear operator defined by
\begin{equation}
  \label{opd}
  P^*_N:\xi_1^{\otimes N}\otimes \xi_2^{\otimes N}\mapsto \Theta_{p;\xi}^{(N)}.
\end{equation}
  Indeed, $P^*_N$  
is the unique bounded linear operator $(L_1^N \boxtimes L_2^N)\Bigr | _{p}\to H^0(M_1\times M_2,L_1^N \boxtimes L_2^N)$ 
determined by the equality
$$
\langle P_N s,v\rangle=\langle s,P_N^*v\rangle
$$
for all $v\in (L_1^N \boxtimes L_2^N)_p$ and all $s\in  H^0(M_1\times M_2,L_1^N \boxtimes L_2^N)$. Since $P_Ns=s(p)$, this equality becomes
\begin{equation}
  \label{opereq}
h^{(N)}(s(p),v)=\int_{M_1\times M_2}h^{(N)}(s(z),P^*_N(v)(z))dV(z).
\end{equation}
Since $v=\lambda \xi_1^{\otimes N}\otimes \xi_2^{\otimes N}$ for some $\lambda\in\C$, it follows that 
 (\ref{opereq}) is equivalent to
$$
h^{(N)}(s(p),\xi_1^{\otimes N}\otimes \xi_2^{\otimes N})=\int_{M_1\times M_2}h^{(N)}(s(z),P_N^*(\xi_1^{\otimes N}\otimes\xi_2^{\otimes N} )(z))dV(z)
$$
for all $s\in  H^0(M_1\times M_2,L_1^N \boxtimes L_2^N)$. 
But, with (\ref{opd}), this is true by (\ref{thetaprop}).

We conclude that the operator $P^*_N P_N$ maps  $H^0(M_1\times M_2,L_1^N \boxtimes L_2^N)$ onto the $1$-dimensional linear subspace
spanned by $\Theta_{p;\xi}^{(N)}$.
Choose an orthonormal basis in 
$H^0(M_1\times M_2,L_1^N \boxtimes L_2^N)$ such that
the first basis vector is $e_1=\frac{1}{||\Theta_{p;\xi}^{(N)}||}\Theta_{p;\xi}^{(N)}$.
Applying (\ref{thetaprop}), we get that all basis vectors, except for the first one, are in the kernel of $P_N$. Therefore,
with respect to this basis, 
the operator $\frac{1}{\tr (P_N^*P_N)}P_N^*P_N$ is the orthogonal projection onto the one-dimensional subspace spanned by $e_1$. Hence, $\rho_N$ is a pure state.  
Moreover, $e_1$ is a decomposable vector, because 
$$
\Theta_{p;\xi}^{(N)}=\pi_1^*\Theta_{p_1;d\pi_1(\xi_1)}^{(N)}\otimes \pi_2^*\Theta_{p_2;d\pi_2(\xi_2)}^{(N)}.
$$
It follows that the entanglement entropy of $\rho_N$ is zero.

Proof of $(ii)$. We have three cases:

 \ $\Lambda_1$ is a point and $\dim(\Lambda_2)>0$,
 
  \ $\dim(\Lambda_1)>0$ and $\Lambda_2$ is a point, 

 \ $\dim(\Lambda_1)>0$ and $\dim(\Lambda_2)>0$.

Suppose $\Lambda_1=\{ p\}$, where $p\in M_1$, and $\dim(\Lambda_2)>0$. Choose and fix a unit vector $\xi\in (L_1)_p$.
Noting (\ref{thetaprop12}), choose an orthonormal basis $e_1$,...,$e_{d_N^{(1)}}$ in $H^0(M_1,L_1^N)$ such that
$$
e_1=\frac{1}{||\Theta_{p,\xi}^{(N)}||}\Theta_{p,\xi}^{(N)}
$$
and $e_j(p)=0$ for $j\ge 2$. Choose an orthonormal basis $f_1$,...,$f_{d_N^{(2)}}$ in  $H^0(M_2,L_2^N)$.
The vectors $e_i\otimes f_j$ form an orthonormal basis in  $H^0(M_1,L_1^N)\otimes H^0(M_2,L_2^N)$
(we choose the following ordering: $e_1\otimes f_1$, $e_1\otimes f_2$,...,$e_1\otimes f_{d_N^{(2)}}$,
$e_2\otimes f_1$, $e_2\otimes f_2$,...,$e_2\otimes f_{d_N^{(2)}}$,...,
$e_{d_N^{(1)}}\otimes f_1$, $e_{d_N^{(1)}}\otimes f_2$,...,$e_{d_N^{(1)}}\otimes f_{d_N^{(2)}}$).
We would like to find the matrix of $\rho_N$ in this basis. 
Because of (\ref{thetaprop12}), $e_i(p)=0$ for $i\ge 2$. Hence, for $i\ge 2$, $P_N(e_i\otimes f_j)=0$, and $\rho_N(e_i\otimes f_j)=0$. 
It remains to determine  $\rho_N(e_1\otimes f_j)$ for each basis vector $f_j$. 
Denote by $d\mu_2$ the measure induced on $\Lambda_2$ by the metric on $M_2$.   
The operator $P_N^*$ is the unique bounded linear operator 
 $$
 V_N\to 
 H^0(M_1\times M_2,L_1^N \boxtimes L_2^N) \cong H^0(M_1,L_1^N)\otimes H^0(M_2,L_2^N)
 $$
 with the property 
 $$
 \langle P_N(e_i\otimes f_j),\xi ^{\otimes N}\otimes f_r|_{\Lambda_2}\rangle=
 \langle e_i\otimes f_j,P_N^*(\xi ^{\otimes N}\otimes f_r|_{\Lambda_2})\rangle
$$
for all $i,j,r$. Since $P_N(e_i\otimes f_j)=e_i(p)\otimes f_j\Bigr|_{\Lambda_2}$, this equality becomes 
$$
\int_{\Lambda_2}h_2^{(N)}(f_j\Bigr| _{\Lambda_2} (z),f_r\Bigr| _{\Lambda_2} (z))d\mu_2 \ h_1^{(N)}(e_i(p),\xi ^{\otimes N})=
$$
$$
\int_{M_1\times M_2}h^{(N)}(e_i(z)\otimes f_j(w),(P_N^*(\xi^{\otimes N}\otimes f_r\Bigr| _{\Lambda_2}))(z,w))dV_1(z)dV_2(w).
$$
Therefore
$$
P_N^*(\xi^{\otimes N}\otimes f_r\Bigr| _{\Lambda_2})=\sum_j \beta_{j}^{(r)}e_1\otimes f_j,
$$
where
$$
\overline{\beta_{j}^{(r)}}=h_1^{(N)}(e_1(p),\xi^{\otimes N})\int_{\Lambda_2}h_2^{(N)}(f_j\Bigr| _{\Lambda_2}(z),f_r\Bigr| _{\Lambda_2}(z))d\mu_2.
$$
Then, since $P_N(e_1)=e_1(p)=c\xi^{\otimes N}$ for some $c\in\C$, we get: 
$$
P_N^*P_N(e_1\otimes f_r)=\sum_j b_{j}^{(r)}e_1\otimes f_j,
$$
where
$$
\overline{b_{j}^{(r)}}=h_1^{(N)}(e_1(p),e_1(p))\int_{\Lambda_2}h_2^{(N)}(f_j\Bigr| _{\Lambda_2}(z),f_r\Bigr| _{\Lambda_2}(z))d\mu_2.
$$
Denote by $B_N$ the $d_N^{(2)}\times d_N^{(2)}$
matrix whose $jr$-th entry 
is $b_{j}^{(r)}$. The orthogonal projector $\Po _{e_1}\in {\mathrm{End}} (H^0(M_1,L_1^N))$, in the basis $\{ e_i\}$, has the matrix which we will denote by 
$E_{11}^{(N)}$, that has $1$ in the $(1,1)$-th entry and zeros in all other entries.  

 The operator $P_N^*P_N$ acts on decomposable vectors by 
$$
a\otimes b\mapsto E_{11}^{(N)}a\otimes B_Nb
$$
for $a\in H^0(M_1,L_1^N)$, $b\in H^0(M_2,L_2^N)$. 
The matrix of $P_N^*P_N$ with respect to the chosen basis is $E_{11}^{(N)}\otimes B_N$ and 
$$
\tr  (P_N^*P_N)=\tr B_N. 
$$
The operator $\rho_N$ acts on decomposable vectors by 
$$
a\otimes b\mapsto E_{11}^{(N)} a\otimes \frac{1}{\tr {B_N}}B_N b
$$
for $a\in H^0(M_1,L_1^N)$, $b\in H^0(M_2,L_2^N)$. 
Each of the matrices $E_{11}^{(N)}$, $\frac{1}{\tr {B_N}}B_N$ is a positive semidefinite Hermitian matrix of trace $1$.   
  Therefore the entanglement of formation of $\rho_N$ is zero. 
 This finishes the proof of
  the first case. The second case is proved similarly. 

The remaining case is  $\dim(\Lambda_1)>0$ and $\dim(\Lambda_2)>0$. 
Choose an orthonormal basis $e_1$,...,$e_{d_N^{(1)}}$ in  $H^0(M_1,L_1^N)$ and an orthonormal basis
$f_1$,...,$f_{d_N^{(2)}}$ in  $H^0(M_2,L_2^N)$.
Then the vectors $e_i\otimes f_j$ form an orthonormal basis in  $H^0(M_1,L_1^N)\otimes H^0(M_2,L_2^N)$ 
(we use the following ordering: $e_1\otimes f_1$, $e_1\otimes f_2$,...,$e_1\otimes f_{d_N^{(2)}}$,
$e_2\otimes f_1$, $e_2\otimes f_2$,...,$e_2\otimes f_{d_N^{(2)}}$,...,
$e_{d_N^{(1)}}\otimes f_1$, $e_{d_N^{(1)}}\otimes f_2$,...,$e_{d_N^{(1)}}\otimes f_{d_N^{(2)}}$). 
We would like to find the matrix of $\rho_N$ in this basis.
For $j=1,2$ denote by $d\mu_j$ the measure induced on $\Lambda_j$ by the metric on $M_j$.
The operator $P_N^*$ is the unique bounded linear operator 
 $$
 V_N\to 
 H^0(M_1\times M_2,L_1^N \boxtimes L_2^N)\simeq H^0(M_1,L_1^N)\otimes H^0(M_2,L_2^N)
 $$
 with the property 
 $$
 \langle P_N(e_i\otimes f_j),e_l|_{\Lambda_1}\otimes f_r|_{\Lambda_2}\rangle=
 \langle e_i\otimes f_j,P_N^*(e_l|_{\Lambda_1}\otimes f_r|_{\Lambda_2})\rangle
$$
for all $i,j,l,r$. Since $P_N(e_i\otimes f_j)=e_i\Bigr |_{\Lambda_1}\otimes f_j\Bigr|_{\Lambda_2}$, this equality becomes 
$$
\int_{\Lambda_1}h_1^{(N)}(e_i\Bigr| _{\Lambda_1} (z),e_l\Bigr| _{\Lambda_1} (z)d\mu_1
\int_{\Lambda_2}h_2^{(N)}(f_j\Bigr| _{\Lambda_2} (z),f_r\Bigr| _{\Lambda_2} (z)d\mu_2=
$$
$$
\int_{M_1\times M_2}h^{(N)}(e_i(z)\otimes f_j(w),(P_N^*(e_l\Bigr| _{\Lambda_1}\otimes f_r\Bigr|_{\Lambda_2}))
(z,w))dV_1(z)dV_2(w).
$$
Hence
$$
P_N^*P_N(e_l\otimes f_r)=\sum_{i,j} a_{i}^{(l)}b_{j}^{(r)}e_i\otimes f_j,
$$
where
$$
\overline{a_{i}^{(l)}}=
\int_{\Lambda_1}h_1^{(N)}(e_i\Bigr| _{\Lambda_1}(z),e_l\Bigr| _{\Lambda_1}(z))d\mu_1,
$$
$$
\overline{b_{j}^{(r)}}=
\int_{\Lambda_2}h_2^{(N)}(f_j\Bigr| _{\Lambda_2}(z),f_r\Bigr| _{\Lambda_2}(z))d\mu_2.
$$
Denote by $A_N$ the $d_N^{(1)}\times d_N^{(1)}$
matrix whose $ij$-th entry 
is $a_{i}^{(j)}$ and denote by $B_N$ the $d_N^{(2)}\times d_N^{(2)}$
matrix whose $ij$-th entry 
is $b_{i}^{(j)}$. The operator $P_N^*P_N$ acts on decomposable vectors by 
$$
a\otimes b\mapsto A_Na\otimes B_Nb
$$
for $a\in H^0(M_1,L_1^N)$, $b\in H^0(M_2,L_2^N)$. 
The matrix of $P_N^*P_N$ with respect to the chosen basis is $A_N\otimes B_N$ and 
$$
\tr  (P_N^*P_N)=(\tr A_N)(\tr B_N). 
$$
The operator $\rho_N$ acts on decomposable vectors by 
$$
a\otimes b\mapsto \frac{1}{\tr {A_N}}A_N \ a\otimes \frac{1}{\tr {B_N}}B_N \ b
$$
for $a\in H^0(M_1,L_1^N)$, $b\in H^0(M_2,L_2^N)$. 
Each of the matrices $\frac{1}{\tr {A_N}}{A_N}$, $\frac{1}{\tr {B_N}}B_N$ is a positive semidefinite Hermitian matrix of trace $1$.   
  Therefore the entanglement of formation of $\rho_N$ is zero. This finishes the proof.  
  $\Box$
  
  Theorem \ref{mainth} asserts that if $\Lambda$ is a product submanifold of $M_1\times M_2$, then, for sufficiently large $N$, the entanglement of formation
  of the associated state $\rho_N$ is zero. The reason for this, as shown in the proof, is that  
  $\rho_N$ is separable for sufficiently large $N$ (i.e. if $\Lambda$ is a product submanifold, then for some $N_0$, the states $\rho_N$ are product states 
  for all $N> N_0$). One could ask whether zero entanglement measure characterizes product submanifolds, or what would characterize product submanifolds. 
  It is not clear how to formulate an ``if and and only if version" of Theorem \ref{mainth}. Indeed, let us consider the following example. Let $\Lambda$ be
  $(M_1\times M_2)-\{ (p_1,p_2)\}$, where $p_1\in M_1$ and $p_2\in M_2$. This submanifold $\Lambda$ is not of the form
  $\Lambda_1\times \Lambda_2$ with $\Lambda_1\subset M_1$ and  $\Lambda_2\subset M_2$.
  However, for sufficiently large $N$, $\tilde{\nu}(\rho_N)=0$, by an argument similar to the one used in the proof of the theorem. Let us explain why. 
The operator $P_N^*$ is the unique bounded linear operator 
 $$
 V_N\to 
 H^0(M_1\times M_2,L_1^N \boxtimes L_2^N)
 $$
 with the property 
 $$
 \langle P_N(s),\tau \Bigr |_{\Lambda}\rangle=
 \langle s,P_N^*(\tau \Bigr |_{\Lambda})\rangle
$$
for all $s,\tau \in H^0(M_1\times M_2,L_1^N \boxtimes L_2^N)$. Since $P_N(s)=s\Bigr |_{\Lambda}$, this equality becomes 
$$
\int_{\Lambda}h^{(N)}(s(z,w),\tau (z,w))dV(z,w)=
\int_{M_1\times M_2}h^{(N)}(s(z,w),(P_N^*(\tau \Bigr |_{\Lambda}))
(z,w))dV_1(z)dV_2(w).
$$
Since $\Lambda$ is $M_1\times M_2$ with one point removed, in the left-hand side
$\int_{\Lambda}=\int_{M_1\times M_2}$. We conclude that 
$P_N^*P_N$ is the identity operator and 
$$
\rho_N=\frac{1}{d_N^{(1)}}I_{1}\otimes \frac{1}{d_N^{(2)}}I_{2},
$$
where $I_j$ is the identity operator on $H^0(M_j,L_j^N)$ for $j=1,2$. 
Therefore $\rho_N$ is separable and hence $\tilde{\nu}(\rho_N)=0$.

\end{document}